\documentclass[12pt,oneside,reqno]{amsart}
\usepackage{txfonts}
\usepackage{bbm}
\usepackage{cases}
\usepackage{amsmath}
\usepackage{graphicx}
\usepackage{mathrsfs}
\usepackage{stmaryrd}
\usepackage{amsfonts}
\usepackage{enumerate,amsmath,amssymb,amsthm}
\numberwithin{equation}{section}

\newcommand{\be}{\begin{eqnarray}}
\newcommand{\ee}{\end{eqnarray}}
\newcommand{\ce}{\begin{eqnarray*}}
\newcommand{\de}{\end{eqnarray*}}
\newtheorem{theorem}{Theorem}[section]
\newtheorem{lemma}[theorem]{Lemma}
\newtheorem{remark}[theorem]{Remark}
\newtheorem{definition}[theorem]{Definition}
\newtheorem{proposition}[theorem]{Proposition}
\newtheorem{Examples}[theorem]{Example}
\newtheorem{corollary}[theorem]{Corollary}

\def\a{\alpha}

\def\p{\partial}

\def\g{\gamma}

\def\[{{\Big[}}
\def\]{{\Big]}}
\def\<{{\langle}}
\def\>{{\rangle}}
\def\({{\Big(}}
\def\){{\Big)}}

\def\bx{{\mathbf{x}}}

\def\dif{{\mathord{{\rm d}}}}
\def\dis{{\mathord{{\rm \bf d}}}}

\def\no{\nonumber}
\def\={&\!\!=\!\!&}
\def\bt{\begin{theorem}}
\def\et{\end{theorem}}
\def\bl{\begin{lemma}}
\def\el{\end{lemma}}
\def\br{\begin{remark}}
\def\er{\end{remark}}
\def\bd{\begin{definition}}
\def\ed{\end{definition}}
\def\bp{\begin{proposition}}
\def\ep{\end{proposition}}
\def\bc{\begin{corollary}}
\def\ec{\end{corollary}}
\def\bx{\begin{Examples}}
\def\ex{\end{Examples}}

\def\cF{{\mathcal F}}

\def\cM{{\mathcal M}}

\def\cT{{\mathcal T}}

\def\mE{{\mathbb E}}

\def\mH{{\mathbb H}}

\def\mN{{\mathbb N}}

\def\mR{{\mathbb R}}

\def\mT{{\mathbb T}}

\def\sC{{\mathscr C}}

\def\fg{{\mathfrak g}}

\def\geq{\geqslant}
\def\leq{\leqslant}

\def\div{\mathord{{\rm div}}}

\allowdisplaybreaks

\begin{document}

\title{Quasi Invariant Stochastic Flows of SDEs with Non-smooth Drifts on Riemannian Manifolds$^*$}

\date{}

\author{Xicheng Zhang}

\thanks{$*$ This work is supported by NSFs of China (No. 10971076; 10871215).}
%Galerkin's Approximation, Analytic Semigroup,

\dedicatory{
School of Mathematics and Statistics, Wuhan University,
Wuhan, Hubei 430072, P.R.China\\
Email: XichengZhang@gmail.com
 }

\begin{abstract}
In this article we prove that stochastic differential equation (SDE)
with Sobolev drift on compact Riemannian manifold admits a unique $\nu$-almost everywhere
stochastic invertible flow, where $\nu$ is the Riemannian measure, which is quasi-invariant with respect to $\nu$. 
In particular, we extend the well known DiPerna-Lions flows of ODEs to SDEs on Riemannian manifold.
\end{abstract}

\maketitle
\rm

\section{Introduction}

Let $M$ be a connected and compact $C^\infty$-manifold of dimension $d$.
Consider the following Stratonovich's stochastic differential equation (SDE) on $M$:
\be
\dif x_t=X_0(x_t)\dif t+X_k(x_t)\circ\dif W^k_t,\ \ x_0=x,\label{SSDE}
\ee
where $X_i, i=0,\cdots,m$ are $m+1$-vector fields on $M$, and $(W_t)_{t\geq 0}$
is the $m$-dimensional standard Brownian motion on the classical Wiener
space $(\Omega,\cF,P; (\cF_t)_{t\geq 0})$, i.e., $\Omega$ is the space of all
continuous functions from $\mR_+$ to $\mR^m$ with locally uniform convergence topology,
$\cF$ is the Borel $\sigma$-field, $P$ is the Wiener measure, $(\cF_t)_{t\geq 0}$
is the natural filtration generated by the coordinate process $W_t(\omega)=\omega(t)$.
Here and below, we use the following convention: if an index appears twice in a product,
it will be summed.

For solving SDE (\ref{SSDE}), there are usually two ways:
One way is to first construct the solutions in local coordinates and then
patches up them (cf. \cite{iw}). Another way is that one embeds $M$ into some Euclidean
space, obtains a solution in this larger space, and then
proves that the solution will actually stay in $M$ if the
starting point $x$ is in $M$ (cf. \cite{El}). Both of these arguments require that
$X_k$, $k=0,\cdots, m$ are smooth (at least $C^2$) vector fields.

In the case of flat Euclidean space, a celebrated theory established by DiPerna and Lions
\cite{Di-Li} says that when $X_0$ only has Sobolev regularity and bounded divergence, ODE
$$
\dif x_t=X_0(x_t)\dif t, x_0=x
$$
defines a unique regular Lagrangian flow in the sense of Lebesgue measure. Their proofs are
based on a new notion called renormalized solution for the associated transport equation:
$$
\p_t u+X_0u=\p_t u+X^i_0\p_iu=0,\ \ u|_{t=0}=u_0,
$$
where $X^i_0$ is the component of vector field $X_0$ under natural frames.
For the DiPerna-Lions flow on compact Riemannian manifold, Dumas, Golse and Lochak
\cite{Du-Go-Lo} gave an outline for the proof.

Recently, we have extended DiPerna-Lions' flow to the case of SDEs in \cite{Zh2}. Therein, we
followed the direct argument of Crippa and De Lellis \cite{Cr-De-Le}. It is worth pointing
out that we can not use the original method of DiPerna and Lions to study the SDEs
with Sobolev drifts
because the associated stochastic partial differential equation is always degenerate
(cf. \cite{Zh2}). On the other hand, when we consider the corresponding SDEs with Sobolev drifts
on Riemannian manifold, it seems that we can not use the localizing and patching method
as well as the embedding method since $X_0$ is not smooth and the solution is only defined for
almost all starting points. In order to extend the result in \cite{Zh2} to
Riemannian manifold, we shall directly use the intrinsic Riemannian distance as in \cite{Re-Zh}.
For this aim, we have to make a detailed analysis for the distance function
associated with the Riemannian metric.

This paper is organized as follows: in Section 2, we give the notion of
$\nu$-almost everywhere
stochastic flow of SDE (\ref{SSDE}) and state our main result. In Section 3,
we analyze the distance function on Riemannian manifold and give some necessary
preliminaries. In Section 4, we prove our main result as in \cite{Cr-De-Le} and \cite{Zh2}
by using the Hardy-Littlewood maximal function on Riemannian manifold.

\section{Main Result}

Let $(M,\fg)$ be a connected and compact $C^\infty$ Riemannian manifold of $d$-dimension,
where $\fg$ denotes the Riemannian metric, a symmetric, positively definite, and
second order covariant tensor field on $M$. Let $\nu(\dif x)$ be the Riemannian measure,
and $\nabla$ the Levi-Civita connection associated with $\fg$.
We also use $\nabla$ to denote the gradient operator.
The divergence operator denoted by
$\div$ is the dual operator of $\nabla$ with respect to $\nu$. Let $TM$ be the tangent bundle.
For any $x\in M$,
the length of a vector $X\in T_xM$ is denoted by $|X|_x:=\sqrt{\fg_x(X,X)}$. Letting $\cT$ be a
measurable transformation of $M$, we use $\nu\circ\cT$ to denote the image measure of $\nu$
under $\cT$, i.e., for any nonnegative measurable function $f$,
$$
\int_M f(x)\nu\circ\cT(\dif x)=\int_M f(\cT(x))\nu(\dif x).
$$
By $\nu\circ\cT\ll\nu$, we mean that $\nu\circ\cT$ is absolutely continuous with respect to $\nu$.

We first introduce the following notion of $\nu$-almost
everywhere stochastic (invertible) flows (cf. \cite{Le-Li} \cite{Am} \cite{Zh2}).
\bd\label{Def1}
Let $x_t(\omega,x)$ be an $M$-valued measurable
stochastic field on $\mR_+\times\Omega\times M$.
We say $x_t(x)$ a {\bf $\nu$-almost everywhere stochastic flow} of
(\ref{SSDE}) corresponding to vector fields $(X_k)_{k=0,\cdots,m}$ if
\begin{enumerate}[{\bf(A)}]
\item For  $\nu$-almost all $x\in\mR^d$, $t\mapsto x_t(x)$ is a continuous and
($\cF_t$)-adapted stochastic process and, satisfies that for any $T>0$ and $f\in C^\infty(M)$,
$$
f(x_t(x))=f(x)+\int^t_0X_0f(x_s(x))\dif s+\int^t_0X_kf(x_s(x))\circ\dif W^k_s,\ \ \forall t\geq 0.
$$
\item For any $t\geq 0$ and $P$-almost all $\omega\in\Omega$,
$\nu\circ x_t(\omega,\cdot)\ll\nu$. Moreover, for any $T>0$, there exists a
constant $K_{T,X_0,X_k}>0$ such that for all nonnegative measurable function $f$ on $M$,
\be\label{Den}
\sup_{t\in[0,T]}\mE\int_M f(x_t(x))\nu(\dif x)\leq K_{T,X_0,X_k}
\int_M f(x)\nu(\dif x).
\ee
\end{enumerate}
We say $x_t(x)$ a {\bf $\nu$-almost everywhere stochastic  invertible flow} of (\ref{SSDE})
corresponding to vector fields $(X_k)_{k=0,\cdots,m}$ if in addition to the above
{\bf (A)} and {\bf (B)},

\begin{enumerate}[{\bf (C)}]
\item  For any $t\geq 0$ and $P$-almost all $\omega\in\Omega$, there exists a measurable
inverse $x^{-1}_t(\omega,\cdot)$ of $x_t(\omega,\cdot)$ so that
$\nu\circ x^{-1}_t(\omega,\cdot)=\rho_t(\omega,\cdot)\nu$,
where the density $\rho_t(x)$ is given by
\be
\rho_t(x):=\exp\left\{\int^t_0\div X_0(x_s(x))\dif s
+\int^t_0\div X_k(x_s(x))\circ\dif W^k_s\right\}.\label{Rho}
\ee
\end{enumerate}
\ed

\br
In the above definitions, we have already assumed that all the integrals make sense.
In particular, the above property {\bf (C)} guarantees the quasi invariance of the flow 
transformation $x\mapsto x_t(x)$ with respect to the Riemannian volume.
\er

For $k\in\mN\cup\{\infty\}$, let $C^k(TM)$ be the set of all
$k$-order smooth vector fields on $M$. For $p\geq 1$ and
$X\in C^\infty(TM)$, we define
$$
\|X\|_p:=\left(\int_M |X|^p_x\nu(\dif x)\right)^{1/p}
$$
and
$$
\|X\|_{1,p}:=\|X\|_p+\left(\int_M |\nabla X|^p_x\nu(\dif x)\right)^{1/p}.
$$
Let $L^p(TM)$ and $\mH^p_1(TM)$ be the completions of $C^\infty(TM)$
with respect to $\|\cdot\|_p$ and $\|\cdot\|_{1,p}$ respectively.
We also use $L^\infty(TM)$ to denote the set of all bounded measurable vector fields.

The following two propositions are direct consequences of Definition \ref{Def1}, whose proofs
can be found in \cite{Zh2}.
\bp\label{Pr2}
Assume that SDE (\ref{SSDE}) admits a unique $\nu$-almost everywhere stochastic flow.
Then the following flow property holds: for any $s\geq 0$ and $(P\times\nu)$-almost all
$(\omega,x)\in\Omega\times M$,
$$
x_{t+s}(\omega,x)=x_t(\theta_s\omega,x_s(\omega,x)),\  \forall t\geq 0,
$$
where $\theta_s\omega:=\omega(s+\cdot)-\omega(s)$. Moreover, for any bounded measurable
function $f$ on $M$, define
$$
\mT_tf(x):=\mE f(x_t(x)),
$$
then for any $t,s\geq 0$
$$
\mE(f(x_{t+s}(x))|\cF_s)=\mT_tf(x_s(x)),\ \ (P\times\nu)-a.e.\label{Ep30}
$$
In particular, $(\mT_t)_{t\geq 0}$ forms a bounded linear operator semigroup
on $L^p(M)$ for any $p\geq 1$.
\ep

\bp 
Assume that $X_0\in L^\infty(TM)$ with $\div X_0\in L^1(M)$
and $X_k\in C^2(TM)$, $k=1,\cdots, m$. Let $x_t(x)$ be a $\nu$-almost everywhere
stochastic invertible flow of (\ref{SSDE}) in the sense of Definition \ref{Def1}.
Let $u_0\in L^1(M)$ and set $u_{t}(x):=u_0(x^{-1}_{t}(x))$.
Then $u_t(x)$ solves the following stochastic transport equation
in the distributional  sense:
$$
\dif u=-X_0 u\dif t-X_k u\circ\dif W^k_t.
$$
In particular,  $\bar u_t(x):=\mE u_0(x^{-1}_t(x))$
is a distributional  solution of the following second order parabolic differential equation:
$$
\p_t \bar u=-\frac{1}{2}\sum_kX_k^2\bar u-X_0 \bar u.
$$
\ep

Our main result in the present paper is:
\bt\label{Main}
Assume that $X_0\in \mH^p_1(TM)\cap L^\infty(TM)$ for some $p>1$, satisfies
$$
\div X_0\in L^\infty(M),
$$
and for each $k=1,\cdots, m$, $X_k\in C^2(TM)$.
Then there exists a unique $\nu$-almost everywhere stochastic invertible
flows $\{x_t(x), x\in M\}_{t\geq 0}$ of  SDE (\ref{SSDE}) in the sense of Definition \ref{Def1}.
\et

\section{Preliminaries}

\subsection{Distance Function}
We  need the following simple lemma.
\bl\label{LL1}
Let $(M,\dis)$ be a compact metric space. Let $\Sigma=\{U_\a,\a\in\Lambda\}$ be a
finite open covering of $M$. Then there exists a positive number $\varrho$ such that
for any $x,y\in M$, if $\dis(x,y)<\varrho$, then $x,y$ must lies in some $U_\a$ simultaneously.
\el
\begin{proof}
We use the contradiction method. Suppose that for any $n\in\mN$, there exists $x_n,y_n\in M$ with
$\dis (x_n,y_n)<\frac{1}{n}$ such that
\be
\mbox{$x_n, y_n$ do not belong to any $U_\a\in\Sigma$ simultaneously}.\label{Con}
\ee
By the compactness of $M$, there is a  subsequence $n_k$ and $z\in M$ such that
$$
\lim_{k\to\infty}\dis (x_{n_k},z)=0,\ \ \lim_{k\to\infty}\dis (y_{n_k},z)=0.
$$
Since $z$ belongs to some open set $U_\a\in\Sigma$, for $k$ large enough, $x_{n_k}$ and
$y_{n_k}$ must lies in $U_\a$, which is contrary to (\ref{Con}). The proof is complete.
\end{proof}

Using this lemma, we have the following property about the distance
function $\dis(\cdot,\cdot)$ on $M$, which will be our localizing basis below.
\bl\label{l1}
Let $M$ be a compact Riemannian manifold. Then, there are a finite covering
$\{(U_\a,\varphi_\a; \xi^k_\a)\}_{\a\in\Lambda}$ of $M$
by normal coordinate neighborhoods, and positive constants $\varrho,\lambda\in(0,1)$ such that
\begin{enumerate}[($1^o$)]
\item For any $x,y\in M$ with $\dis(x,y)<\varrho$, $x,y$ must be in some $U_\a$ simultaneously,
and there is a unique minimizing geodesic connecting $x$ and $y$ in $U_\a$.

\item
In local coordinate $\{(U_\a,\varphi_\a; \xi^k_\a\}$, for any $x,y\in U_\a$,
$$
\lambda\cdot|\varphi_\a(x)-\varphi_\a(y)|\leq\dis(x,y)\leq
\lambda^{-1}\cdot|\varphi_\a(x)-\varphi_\a(y)|,
$$
where $|\cdot|$ denotes the Euclidean metric in $\varphi_\a(U_\a)\subset\mR^d$. Moreover,
$$
\lambda I\leq (g^\a_{ij})\leq \lambda^{-1}I,
$$
where $g^\a_{ij}:=\fg(\partial_{\xi^i_\a}, \partial_{\xi^j_\a})$.

\item For any $U_\a$, the restriction of $\dis^2(\cdot,\cdot)$ to
$U_\a\times U_\a$ belongs to $C^\infty(U_\a\times U_\a)$.
\end{enumerate}
\el
\begin{proof}
For each $a\in M$, there is a normal coordinate neighborhood $(U_a,\varphi_a)$ of $a$
such that any two points in $U_a$ can be joined by a unique minimizing geodesic lying in $U_a$,
and $\dis^2(\cdot,\cdot)\in C^\infty(U_a\times U_a)$ (see \cite[p.166, Theorem 3.6]{K-N}).
Moreover, there is a constant $\lambda_a$ such that
for all $x,y\in U_a$ (see \cite[p.125]{B-C}),
$$
\lambda_a^{-1}|\varphi_a(x)-\varphi_a(y)|\leq \dis(x,y)\leq
\lambda_a|\varphi_a(x)-\varphi_a(y)|
$$
and
$$
\lambda_a^{-1}I\leq (g^\a_{ij})\leq \lambda_aI,
$$
The results now follow by the compactness of $M$ and Lemma \ref{LL1}.
\end{proof}

In the following, we shall fix the $\Sigma:=\{(U_\a,\varphi_\a;
\xi^k_\a)\}_{\a\in\Lambda}$ and $\varrho,\lambda$ in this lemma as well
as a unit partition $(\psi_\a)_{\a\in\Lambda}$
subordinate to $\Sigma$, i.e.,
\be
\psi_\a \in C^\infty(M;[0,1]),\ \ \mbox{supp}(\psi_\a)\subset U_\a,\ \
\sum_{\a\in\Lambda}\psi_\a\equiv 1.\label{Unit}
\ee

Given two points $x,y\in M$ with $\dis(x,y)<\varrho$, let
$$
\{\g(s),\ s\in[0,t_0],\ t_0:=\dis(x,y)\}
$$
with
$$
\g(0)=x,\ \ \g(t_0)=y
$$
be the unique minimizing geodesic connecting $x$ and $y$. We use $//^\g_{y\leftarrow x}$ to
denote the parallel transport from $x$ to $y$ along the geodesic
$\g$, i.e., $//^\g_{y\leftarrow x}$ establishes an isomorphism between
tangent spaces $T_xM$ and $T_{y}M$. For a vector field $X$ and a smooth function $f$, we write
$$
\fg_x(X,\nabla f)=X(x)f=[Xf](x).
$$

\bl\label{l2} For $x\in M$ and a vector $X\in T_xM$, we have
$$
\fg_x(X,\nabla\dis(\cdot,y))=-\fg_y(//^\g_{y\leftarrow x}X,\nabla\dis(x,\cdot)).
$$
\el
\begin{proof}
By a corollary to Gaussian Lemma (see e.g. \cite[Corollary
6.9]{lee}), we have
\ce
\fg_x(X,\nabla\dis(\cdot,y))&=&\fg_x(X,\dot{\g}(0))=-\fg_y(//^\g_{y\leftarrow x}X,\dot{\g}(t_0))\\
&=&-\fg_y(//^\g_{y\leftarrow x}X,\nabla \dis(x,\cdot)).
\de
\end{proof}

\subsection{Local maximal function on Riemannian manifold $M$}

\textsc{Convention}:
For two expressions $A$ and $B$, the notation $A\preceq B$ means that
$A\leq C\cdot B$, where $C>0$ is an unimportant constant and
may change in different occasions. We assume that the reader can see
the dependence of $C$ on the parameters from the context.

For a nonnegative function $f\in L^1(M)$ and $R>0$,
the local maximal function $\cM_Rf$ is defined by
$$
\cM_R f(x):=\sup_{r\in(0,R)}\frac{1}{\nu(B_r(x))}\int_{B_r(x)}f(y)\nu(\dif y),
$$
where $B_r(x):=\{y\in M: \dis (x,y)< r\}$.
Similarly, for a function $h\in L^1_{loc}(\mR^d)$,
we  define the local maximal function $\tilde \cM_Rh$ in Euclidean space $\mR^d$ by
$$
\tilde\cM_R h(\xi):=\sup_{r\in(0,R)}\frac{1}{|\tilde B_r(\xi)|}\int_{\tilde B_r(\xi)}h(\eta)\dif \eta,
$$
where $\tilde B_r(\xi):=\{\eta\in \mR^d: |\eta-\xi|<r\}$ and
$|\tilde B_r(\xi)|$ denotes the volume of ball $\tilde B_r(\xi)$ with respect to
the Lebesgue measure.

We have
\bl\label{Le3}
Let $f$ be a measurable function on $M$ with $\nabla f\in L^1(TM)$. Then,
there exists a $\nu$-null set $N$ such that for
all $x,y\notin N$ with $\dis (x,y)<\lambda^2\varrho$,
$$
|f(x)-f(y)|\preceq\dis (x,y)\cdot(\cM_{\varrho}|\nabla f|(x)+\cM_{\varrho}|\nabla f|(y)),
$$
where $\lambda$ and $\varrho$ are from Lemma \ref{l1}.
\el
\begin{proof}
Since $\dis(x,y)<\varrho$, by ($1^o$) of Lemma \ref{l1} we only need to prove the lemma in
local coordinate $(U,\varphi; \xi^k)\in\Sigma$.
It is well known that there is a Lebesgue-null set $Q$ such that for all
$\xi,\eta\in \varphi(U)\setminus Q$ with $|\xi-\eta|<\lambda\varrho$
(cf. \cite[Appendix]{Cr-De-Le}),
$$
|f\circ\varphi^{-1}(\xi)-f\circ\varphi^{-1}(\eta)|\preceq |\xi-\eta|
\cdot(\tilde\cM_{\lambda\varrho}|\nabla (f\circ\varphi^{-1})|(\xi)
+\tilde\cM_{\lambda\varrho}|\nabla (f\circ\varphi^{-1})|(\eta)).
$$
Noting that by ($2^o$) of Lemma \ref{l1},
$$
\varphi(B_{\lambda r}(x))\subset\tilde B_r(\varphi(x))\subset\varphi(B_{\lambda^{-1} r}(x))
$$
and
$$
\lambda^{d/2}\nu(B_{\lambda r}(x))\leq|\tilde B_r(\varphi(x))|\leq
\lambda^{-d/2}\nu(B_{\lambda^{-1} r}(x)),
$$
we thus have
$$
\tilde\cM_{\lambda\varrho}|\nabla (f\circ\varphi^{-1})|(\xi)
\preceq\cM_{\varrho}|\nabla f|(\varphi^{-1}(\xi)).
$$
The result now follows.
\end{proof}

The following result can be proved along the same lines as in \cite[p.5 Theorem 1]{St}.
\bl
Let $f\in L^p(M)$ for some $p>1$, then
\be
\|\cM_R f\|_p\preceq \|f\|_p.\label{El1}
\ee
\el
\subsection{Two estimates about vector fields}

\bl\label{Le0}
Let $X\in \mH^1_1(TM)$ be a Sobolev vector field.
Then there exists a $\nu$-null set $N$ such that for all $x,y\notin N$
with $\dis(x,y)<\lambda^2\varrho$,
\be
|X(x)\dis^2(\cdot,y)+X(y)\dis^2(x,\cdot)|\preceq \dis^2(x,y)\cdot
(1+\cM_\varrho|X|_1(x)+\cM_\varrho|X|_1(y)),
\ee
where $|X|_1(x):=|X|_x+|\nabla X|_x$, and the constant in $\preceq$ is independent of $X$.
In particular, if $X\in C^1(TM)$, then
\be
|X(x)\dis^2(\cdot,y)+X(y)\dis^2(x,\cdot)|\preceq \dis^2(x,y)
\cdot\Big(2\sup_{x\in M}|X|_1(x)+1\Big).\label{Ep3}
\ee
\el
\begin{proof}
By Lemma \ref{l2}, we have
\ce
|X(x)\dis^2(\cdot,y)+X(y)\dis^2(x,\cdot)|
&=&2\dis(x,y)\cdot|\fg_y(X(y)-//^\g_{y\leftarrow x}X(x),\nabla\dis(x,\cdot))|\\
&\leq&2\dis(x,y)\cdot|X(y)-//^\g_{y\leftarrow x}X(x)|_{y}.
\de
Thus, it is enough to prove that there exists a $\nu$-null set $N$ such that for all $x,y\notin N$
with $\dis(x,y)<\lambda^2\varrho$,
$$
|X(y)-//^\g_{y\leftarrow x}X(x)|_{y}\preceq \dis(x,y)\cdot
(1+\cM_\varrho|X|_1(x)+\cM_\varrho|X|_1(y)).
$$
Since $\dis(x,y)<\lambda^2\varrho$, we only need to prove it in a
local coordinate $(U,\varphi; \xi^k)\in\Sigma$. In  local coordinate
$(U,\varphi; \xi^k)$, we may write
$$
X(x)|_U=X^k(x)\p_{\xi^k}
$$
and
$$
\nabla X(x)|_U=(\p_{\xi^i} X^k+X^j\Gamma^k_{ji})\dif \xi^i\otimes\p_{\xi^k},
$$
where $\Gamma^k_{ji}=\fg(\nabla_{\p_{\xi^j}}\p_{\xi^i},\p_{\xi^k})$
are Christoffel symbols.
By Lemma \ref{Le3}, there exists a $\nu$-null set $N$ such that for all $x,y\notin U\setminus N$ with
$\dis (x,y)<\lambda^2\varrho$,
\be
|X^k(x)-X^k(y)|&\preceq& \dis(x,y)\cdot(\cM_\varrho|\nabla X^k|(x)
+\cM_\varrho|\nabla X^k|(y))\no\\
&\preceq&\dis(x,y)\cdot(\cM_\varrho|X|_1(x)+\cM_\varrho|X|_1(y)).\label{Ep1}
\ee

Let $t_0:=\dis(x,y)$, and $\{Y^k_s,s\in[0,t_0],k=1,\cdots,d\}$ be the unique solution to ODEs
$$
\frac{\dif Y^k_s}{\dif s}+\sum_{ij}\Gamma^k_{ij}(\g(s))\cdot Y^i_s\cdot \dot{\g}^j_s=0,\
\ Y^k_0=X^k(x), \ k=1,\cdots,d.
$$
Then $//^\g_{y\leftarrow x}X(x)=Y^k_{t_0}\cdot\p_{\xi^k}$. From this equation, one easily finds that
\be
|Y^k_{t_0}-X^k(x)|=|Y^k_{t_0}-Y^k_{0}|\preceq t_0=\dis(x,y).\label{Ep2}
\ee
Hence,  by ($2^o$) of Lemma \ref{l1},
\ce
|X(y)-//^\g_{y\leftarrow x}X(x)|_{y}&=&\left((X^k(y)-Y^k_{t_0})
\cdot(X^j(y)-Y^j_{t_0})\cdot g_{kj}(y)\right)^{1/2}\\
&\preceq&\sum_{k=1}^d|X^k(y)-Y^k_{t_0}|\\
&\preceq& \sum_{k=1}^d\left(|X^k(y)-X^k(x)|+|X^k(x)-Y^k_{t_0}|\right)\\
&\preceq& \dis(x,y)\cdot(1+\cM_\varrho|X|_1(x)+\cM_\varrho|X|_1(y)),
\de
where the last step is due to (\ref{Ep1}) and (\ref{Ep2}). The proof is finished.
\end{proof}
\bl\label{l3} Let $X$ be a $C^2$-vector field on $M$.
Then for any $x,y\in M$ with $\dis(x,y)<\varrho$,
\be
|(X^2\dis^2)_{11}(x,y)+(X^2\dis^2)_{12}(x,y)+(X^2\dis^2)_{21}(x,y)
+(X^2\dis^2)_{22}(x,y)|\preceq\dis^2(x,y),\label{e2}
\ee
where $(X^2\dis^2)_{12}(x,y)=X(y)X(x)\dis^2(x,y)$ and similarly for others, and the constant
in $\preceq$ may depend on $X$.
\el
\begin{proof}
First of all, we have
$$
(X^2\dis^2)_{11}(x,y)=X\fg(X,\nabla\dis^2(\cdot,y))(x)
=\fg_x(\nabla_{X}X,\nabla\dis^2(\cdot,y))+\fg_x(X,\nabla_{X}\nabla\dis^2(\cdot,y))
$$
and
$$(X^2\dis^2)_{22}(x,y)=X\fg(X,\nabla\dis^2(x,\cdot))(y)
=\fg_y(\nabla_{X}X,\nabla\dis^2(x,\cdot))+\fg_y(X,\nabla_{X}\nabla\dis^2(x,\cdot)),
$$
where the second equality is due to the property of the Levi-Civia connection.

By Lemma \ref{l2}, we also have
\ce
(X^2\dis^2)_{12}(x,y)&=&X\fg_x(X,\nabla\dis^2(\cdot,y))(y)
=-X\fg(//^\g_{\cdot\leftarrow x}X,\nabla\dis^2(x,\cdot))(y)\\
&=&-\fg_y(\nabla_{X}(//^\g_{\cdot\leftarrow x}X),\nabla\dis^2(x,\cdot))-
\fg_y(//^\g_{\cdot\leftarrow x}X,\nabla_{X}\nabla\dis^2(x,\cdot))
\de
and
\ce
(X^2\dis^2)_{21}(x,y)&=&X\fg_y(X,\nabla\dis^2(x,\cdot))(x)
=-X\fg(//^\g_{\cdot\leftarrow y}X,\nabla\dis^2(\cdot,y))(x)\\
&=&-\fg_x(\nabla_{X}(//^\g_{\cdot\leftarrow y}X),\nabla\dis^2(\cdot,y))-
\fg_x(//^\g_{\cdot\leftarrow y}X,\nabla_{X}\nabla\dis^2(\cdot,y)).
\de
Thus,
\ce
(X^2\dis^2)_{11}(x,y)+(X^2\dis^2)_{12}(x,y)+(X^2\dis^2)_{21}(x,y)
+(X^2\dis^2)_{22}(x,y)=\mbox{I}+\mbox{II}+\mbox{III},
\de
where
\ce
\mbox{I}&:=&\fg_x(\nabla_{X} X,\nabla\dis^2(\cdot,y))+\fg_y(\nabla_{X}X,\nabla\dis^2(x,\cdot))\\
&=&\fg_x(\nabla_{X} X-//^\g_{\cdot\leftarrow y}\nabla_{X}X,\nabla\dis^2(\cdot,y)),\\
\mbox{II}&:=&-\fg_x(\nabla_{X}(//^\g_{\cdot\leftarrow y}X),
\nabla\dis^2(\cdot,y))-\fg_y(\nabla_{X}(//^\g_{\cdot\leftarrow x}X),\nabla\dis^2(x,\cdot))\\
&=&\fg_x(//^\g_{x\leftarrow y}(\nabla_{X}(//^\g_{\cdot\leftarrow x}X))
-\nabla_{X}(//^\g_{\cdot\leftarrow y}X), \nabla\dis^2(\cdot,y)),\\
\mbox{III}&:=&\fg_x(X,\nabla_{X}\nabla\dis^2(\cdot,y))
+\fg_y(X,\nabla_X\nabla\dis^2(x,\cdot))\\
&&-\fg_x(//^\g_{\cdot\leftarrow y}X,\nabla_{X}\nabla\dis^2(\cdot,y))
-\fg_y(//^\g_{\cdot\leftarrow x}X,\nabla_{X}\nabla\dis^2(x,\cdot))\\
&=&\fg_x(X-//^\g_{x\leftarrow y}X,\nabla_{X}\nabla\dis^2(\cdot,y)
-//^\g_{x\leftarrow y}\nabla_{X}\nabla\dis^2(x,\cdot)).
\de
Now, in a local coordinate $(U,\varphi; \xi^k)\in\Sigma$, set for $k=1,\cdots, d$
\ce
Z^k_1(x,y)&:=&(\nabla_{X}(//^\g_{\cdot\leftarrow x}X(x)))^k(y),\\
Z^k_2(x,y)&:=&(\nabla_{X}\nabla\dis^2(\cdot,y))^k(x).
\de
It is easy to see that $Z^k_1$ and $Z^k_2$ are $C^1$ functions on $U\times U$. Hence,
\ce
|Z^k_1(x,y)-Z^k_1(y,x)|&\preceq&  \dis(x,y),\\
|Z^k_2(x,y)-Z^k_2(y,x)|&\preceq&  \dis(x,y).
\de
As in the proof of Lemma \ref{Le0}, one has
$$
|//^\g_{x\leftarrow y}(\nabla_{X}(//^\g_{\cdot\leftarrow x}X))-
\nabla_{X}(//^\g_{\cdot\leftarrow y}X)|_{x}\preceq \dis(x,y)
$$
and
$$
|\nabla_{X}\nabla\dis^2(\cdot,y)
-//^\g_{x\leftarrow y}(\nabla_{X}\nabla\dis^2(x,\cdot))|_{x}\preceq \dis(x,y).
$$
Combining (\ref{Ep3}) and the above estimates, we obtain the desired result.
\end{proof}

\subsection{Mollifying a non-smooth vector field}

For any measurable vector field $X\in TM$, recalling (\ref{Unit}), we may write
$$
X=\sum_\a \psi_\a X=\sum_\a \psi_\a X|_{U_\a}=\sum_{\a}\psi_\a X^k_\a\p_{\xi^k_\a},
$$
where $X^k_\a: U_\a\to\mR$ is the coordinate component of $X$ in local coordinate
$(U_\a,\varphi_\a; \xi^k_\a)$.

Let $\zeta$ be a nonnegative smooth function on $\mR^d$ with support in
$\{\xi\in\mR^d: |\xi|<1\}$ and
$$
\int_{\mR^d}\zeta(\xi)\dif \xi=1.
$$
Set $\zeta_n(\xi):=n^{d}\zeta(n\xi)$ and define
\be
X^{k}_{\a,n}(x):=(X^k_\a\circ\varphi^{-1}_\a*\zeta_n)\circ\varphi_\a
:=\int_{\varphi_\a(U_\a)}X^k_\a\circ\varphi^{-1}_\a(\xi)
\cdot\zeta_n(\varphi_\a(x)-\xi)\dif \xi\label{Mol2}
\ee
and
\be
X_n:=\sum_\a\psi_\a X^{k}_{\a,n}\p_{\xi^k_\a}.\label{Mol}
\ee
Then it is clear that $X_n\in C^\infty(TM)$.
\br
In general, the restriction of $X_n$ to $U_\a$ does not equal to $X^{k}_{\a,n}\p_{\xi^k_\a}$
since for $\a\not=\beta$, the following compatibility is not true any more:
$$
X^{k}_{\a,n}\not=X^{j}_{\beta,n}\p \xi^k_\a/\p \xi^j_\beta
\ \ \mbox{ in  }U_\a\cap U_\beta\not=\emptyset.
$$
\er
We have the following proposition.
\bp\label{Pro2}
Let $X\in\mH^p_1(TM)$ for some $p\geq 1$ and $X_n$ be defined by (\ref{Mol}). Then
$$
\lim_{n\to \infty}\|X-X_n\|_{1,p}=0.
$$
Moreover, if $X\in L^\infty(TM)$ satisfies $[\div X]^-\in L^\infty(M)$,
then for some constant $C>0$ independent of $n$ and $X$,
\be
\|[\div X_n]^-\|_{L^\infty(M)}\leq C(\|[\div X]^-\|_{L^\infty(M)}
+\|X\|_{L^\infty(TM)}).\label{PP1}
\ee
\ep
\begin{proof}
First of all, by ($2^o$) of Lemma \ref{l1}, we have
\ce
\lim_{n\to \infty}\|X-X_n\|^p_{p}&\leq&\lim_{n\to \infty}\sum_\a\int_{U_\a}
\psi^p_\a [(X^k_\a-X^k_{\a,n})(X^j_\a-X^j_{\a,n})g^\a_{kj}]^{p/2}\nu(\dif x)\\
&\leq&C\lim_{n\to \infty}\sum_{\a,k}\int_{\varphi_\a(U_\a)}
|X^k_\a\circ\varphi^{-1}_\a-X^k_\a\circ\varphi^{-1}_\a*\zeta_n|^p\dif \xi=0.
\de
Similarly, one has
$$
\lim_{n\to \infty}\|\nabla(X-X_n)\|^p_{p}=0.
$$
Moreover, noting that
$$
\div X|_{U_\a}=X^{k}_{\a}\Gamma^i_{ki}
+\p_{\xi^\a_k} X^k_{\a},
$$
we have
$$
\|[\p_{\xi^\a_k} X^k_{\a}]^-\|_{L^\infty(U_\a)}\leq\|[\div X]^-\|_{L^\infty(M)}
+C\|X\|_{L^\infty(TM)}.
$$
Thus, by (\ref{Mol2}) we have
\ce
\|[\div X_n]^-\|_{L^\infty(M)}&\leq&\Big\|\sum_\a \left(\psi_\a (X^{k}_{\a,n}\Gamma^i_{ki}
+\p_{\xi^k_\a} X^k_{\a,n})+X^k_{\a,n}\p_{\xi_k^\a}\psi_\a\right)^-\Big\|_{L^\infty(M)}\\
&\leq&\sum_\a \left(\psi_\a (\|X^{k}_{\a,n}\Gamma^i_{ki}\|_{L^\infty(U_\a)}
+\|[\p_{\xi^k_\a} X^k_{\a,n}]^-\|_{L^\infty(U_\a)})
+\|X^k_{\a,n}\p_{\xi^k_\a}\psi_\a\|_{L^\infty(U_\a)}\right)\\
&\leq& C(\|[\div X]^-\|_{L^\infty(M)}+\|X\|_{L^\infty(TM)}).
\de
The proof is complete.
\end{proof}

\section{Proof of Main Result}

We first prove the following key estimation.
\bl\label{Le2}
Let $x_t(x)$ and $\hat x_t(x)$ be two $\nu$-almost everywhere
stochastic flows of (\ref{SSDE}) corresponding to $(X_0, X_k, k=1,\cdots,m)$ and
$(\hat X_0, X_k, k=1,\cdots,m)$, where
$$
X_0,\hat X_0\in\mH^{p}_1(TM)\mbox{ for some $p>1$ and }X_k\in C^2(TM),\ \ k=1,\cdots,m.
$$
Then for any $\delta>0$,
$$
\mE\int_M\log\left(\frac{\sup_{t\in[0,T]}\dis^2(x_t(x),\hat
x_t(x))}{\delta^2}+1\right) \nu(\dif x)\leq C_1+\frac{C_2}{\delta}
\|X_0-\hat X_0\|_1,
$$
where $C_1= C\cdot (1+K_{T,X_0,X_k}+K_{T,\hat X_0,X_k})(1+\|X_0\|_{1,p})$
and $C_2= C\cdot K_{T,\hat X_0,X_k}$.
Here, $K_{T,X_0,X_k}$ is from (\ref{Den}) and the constant $C$ is independent of
$\delta$ and $X_0, \hat X_0$.
\el
\begin{proof}
Below, let $\chi:\mR^+\to\mR^+$ be a smooth function satisfying
$$
\chi(r)=r, \ \ r\in[0,\lambda^4\varrho^2/4];\ \ \chi(r)=\lambda^4\varrho^2/2,\ \
r\in[\lambda^4\varrho^2,\infty).
$$
We define
$$
f(x,y):=\chi(\dis^2(x,y)).
$$
Then by Lemma \ref{l1}, $f\in C^\infty(M\times M)$ satisfies
$$
f(x,y)\leq\dis^2 (x,y)\leq C_{\varrho,\lambda} f(x,y).
$$
For the simplicity of notations, we write $z_t(x):=(x_t(x),\hat x_t(x))$.
By It\^o's formula, we have
\ce
f(z_t(x))&=&\int^t_0[(X_0f)_1 +(\hat
X_0f)_2](z_s(x))\dif s
+\int^t_0[(X_kf)_1+(X_kf)_2](z_s(x))\circ\dif W^k_s\\
&=&\int^t_0[(X_0f)_1+(\hat X_0f)_2](z_s(x))\dif s
+\int^t_0[(X_kf)_1+(X_kf)_2](z_s(x))\dif W^k_s\\
&&+\frac{1}{2}\int^t_0[(X^2_kf)_{11}+(X_k^2f)_{21} +(X_k^2f)_{12}+(X^2_kf)_{22}](z_s(x))\dif s,
\de
where $(X_0 f)_1(x,y)=X_0(x)f(\cdot,y)$ and similarly for others.
Using It\^o's formula again, we further have
\ce
\log\left(\frac{f(z_t(x))}{\delta^2}+1\right)&=&\int^t_0\frac{[(X_0f)_1
+(\hat X_0f)_2](z_s(x))}{f(z_s(x))+\delta^2}\dif s
+\int^t_0\frac{[(X_kf)_1+(X_kf)_2](z_s(x))}
{f(z_s(x))+\delta^2}\dif W^k_s\\
&&+\frac{1}{2}\int^t_0\frac{[(X^2_kf)_{11}+(X_k^2f)_{21} +(X_k^2f)_{12}+(X^2_kf)_{22}](z_s(x))}{f(z_s(x))+\delta^2}\dif s\\
&&-\frac{1}{2}\int^t_0\frac{|[(X_kf)_1+(X_kf)_2](z_s(x))|^2}
{(f(z_s(x))+\delta^2)^2}\dif s\\
&=:&I_1(t,x)+I_2(t,x)+I_3(t,x)+I_4(t,x).
\de

Let us first treat $I_1(t,x)$. We write
\ce
I_1(t,x)&=&\int^t_0\frac{[(X_0f)_1
+(X_0f)_2](z_s(x))}{f(z_s(x))+\delta^2}\dif s
+\int^t_0\frac{[(\hat X_0f)_2-(X_0f)_2](z_s(x))}
{f(z_s(x))+\delta^2}\dif s\\
&=:&I_{11}(t,x)+I_{12}(t,x).
\de
For a continuous real function $h(t)$, we write
$$
h^*(T):=\sup_{t\in[0,T]}h(t).
$$
By Lemma \ref{Le0}, we have
\ce
\mE\int_M I^*_{11}(T,x)\nu(\dif x)&\preceq&
\mE\int_M\!\int^T_0\frac{\chi'(\dis^2(z_s(x)))\cdot|[(X_0\dis^2)_1
+(X_0\dis^2)_2](z_s(x))|}{\dis^2(z_s(x))+\delta^2}\dif s\nu(\dif x)\\
&\preceq&\mE\int^T_0\!\!\!\int_M(1+\cM_\varrho|X_0|_1(x_s(x))
+\cM_\varrho|X_0|_1(\hat x_s(x)))\nu(\dif x)\dif s\\
&\stackrel{(\ref{Den})}{\preceq}&
(K_{T,X_0,X_k}+K_{T,\hat X_0,X_k})\int_M(1+\cM_\varrho|X_0|_1(x))\nu(\dif x)\\
&\stackrel{(\ref{El1})}{\preceq}&
(K_{T,X_0,X_k}+K_{T,\hat X_0,X_k})(1+\|X_0\|_{1,p}).
\de
Noticing that
\ce
|(\hat X_0f)_2-(X_0f)_2|(x,y)&=&|\chi'(\dis^2(x,y))\cdot((\hat X_0 \dis^2)_2-
(X_0 \dis^2)_2)(x,y)|\\
&\leq&|\chi'(\dis^2(x,y))|\cdot\dis(x,y)\cdot|\hat X_0(y)-X_0(y)|_y,
\de
we similarly have
\ce
\mE\int_M I^*_{12}(T,x)\nu(\dif x)&\preceq&
\frac{1}{\delta}\mE\int^T_0\!\!\!\int_M|\hat X_0(\hat x_s(x))-X_0(\hat x_s(x))|_{\hat x_s(x)}\nu(\dif x)\dif s\\
&\preceq& \frac{K_{T,\hat X_0, X_k}}{\delta}\int_M|X_0-\hat X_0|_x\nu(\dif x).
\de
For $I_2(t,x)$, by BDG's inequality and Lemma \ref{Le0}, we have
$$
\mE\int_M I^*_{2}(T,x)\nu(\dif x)\preceq\int_M\mE\left(\int^T_0
\left|\frac{[(X_kf)_1+(X_kf)_2](z_s(x))}
{f(z_s(x))+\delta^2}\right|^2\dif s\right)^{1/2}\nu(\dif x)\leq C,
$$
where the constant $C$ is independent of $\delta$ and may depend on $X_k$.
Similarly, by Lemma \ref{l3}, we also have
$$
\mE\int_M I^*_{3}(T,x)\nu(\dif x)\leq C.
$$
Since $I_4(t,x)$ is negative, this term can be dropped. Combining the above calculations,
we obtain the desired estimate.
\end{proof}

We also recall the following results for later use (cf. \cite{Zh2}).
\bl\label{Le4}
Let $x_n(\omega,x):\Omega\times M\to M, n\in\mN$
be a family of measurable mappings. Suppose that for $P$-almost all $\omega\in\Omega$,
$\nu\circ x_n(\omega,\cdot)\ll\nu$ and the density $\beta_n(\omega,x)$ satisfies
\be
\sup_n \sup_{x\in M}\mE|\beta_n(x)|^2\leq C_1.\label{Es66}
\ee
If for ($P\times\nu$)-almost all $(\omega,x)\in\Omega\times M$,
$x_n(\omega,x)\to x_0(\omega,x)$ as
$n\to\infty$, then for $P$-almost all $\omega\in\Omega$, $\nu\circ x_0(\omega,\cdot)\ll
\nu$ and the density $\beta$ also satisfies
\be
 \sup_{x\in M}\mE|\beta(x)|^2\leq C_1.\label{Es77}
\ee
Moreover, let $(f_n)_{n\in\mN}$ be a family of uniformly bounded and measurable functions
on $M$. If $f_n$ converges to some $f$ in $L^1(M)$, then
\be
\lim_{n\to\infty}\mE\int_M|f_n(x_n(x))-f(x_0(x))|\nu(\dif x)=0.\label{Ps2}
\ee
\el
\bl\label{Le33}
Let $\cT,\hat\cT:M\to M$ be two measurable transformations. Let $\sC$ be a countable
and dense subset of $C(M)$. Let $\rho\in L^1(M)$ be a positive measurable function.
Assume that for any $f,g\in \sC$,
$$
\int_M f(\hat\cT(x))\cdot g(x)\nu(\dif x)
=\int_M f(x)\cdot g(\cT(x))\cdot\rho(x)\nu(\dif x).\label{L11}
$$
Then $\cT$ admits a measurable invertible $\hat\cT$, i.e., $\cT^{-1}(x)=\hat\cT(x)$ a.e..
Moreover,
$$
\nu\circ \cT^{-1}=\rho\nu,\ \ \ \nu\circ \cT=\rho^{-1}(\cT^{-1})\nu.
$$
\el

\bp\label{Pro1}
Consider SDE (\ref{SSDE}) with $X_k\in C^2(TM)$, $k=0,1,\cdots,m$.
Let $x_t(x)$ be the unique stochastic homeomorphism flow associated with SDE (\ref{SSDE}). Then
$$
\nu\circ x^{-1}_t(\dif x)\sim\nu(\dif x),\ \ \ \nu\circ x_t(\dif x)\sim\nu(\dif x)
$$
and
\be
\nu\circ x^{-1}_t(\dif x)=\exp\left\{\int^t_0\div X_0(x_s(x))\dif s
+\int^t_0\div X_k(x_s(x))\circ\dif W^k_s\right\}\nu(\dif x).\label{Es1}
\ee
Moreover, for any $q\geq 1$
\be
\mE\left|\frac{\nu\circ x^{-1}_t(\dif x)}{\nu(\dif x)}\right|^q
\leq \exp\left\{C_q T(\|[\div X_0]^+\|_\infty
+\|\div X_k\|^2_\infty+\|X_k\div X_k\|_\infty)\right\}\label{Es2}
\ee
and
\be
\mE\left|\frac{\nu\circ x_t(\dif x)}{\nu(\dif x)}\right|^q
\leq \exp\left\{C_q T(\|[\div X_0]^-\|_\infty
+\|\div X_k\|^2_\infty+\|X_k\div X_k\|_\infty)\right\}.\label{Es22}
\ee
\ep
\begin{proof}
We sketch the proof. Let $W_{n,t}$ be the linearized approximation of $W_t$.
Consider the following ODE on $M$:
$$
\dif x_{n,t}(x)=X_0(x_{n,t}(x))\dif t+X_k(x_{n,t}(x)) \dot W_{n,t}^k\dif t.
$$
It is a well known fact that
$$
\nu\circ x^{-1}_{n,t}(\dif x)=\exp\left\{\int^t_0\div X_0(x_{n,s}(x))\dif s
+\int^t_0\div X_k(x_{n,s}(x))\dot W_{n,s}^k\dif s\right\}\nu(\dif x)
$$
By the limit theorem (cf. \cite{Ku}, \cite{Str}, \cite{Re-Zh}),
the desired formula (\ref{Es1}) then follows.

Note that
$$
\int^t_0\div X_k(x_s(x))\circ\dif W^k_s=\int^t_0\div X_k(x_s(x))\dif W^k_s
+\frac{1}{2}\int^t_0 X_k\div X_k(x_s(x))\dif s
$$
and
$$
t\mapsto\exp\left\{q\int^t_0\div X_k(x_s(x))\dif W^k_s-\frac{q^2}{2}
\int^t_0|\div X_k|^2(x_s(x))\dif s\right\}
$$
is an exponential martingale. It is easy to see that (\ref{Es2}) holds.
(\ref{Es22}) can be proved similarly (cf. \cite{Zh2}).
\end{proof}

We now prove the following result.
\bt\label{Me}
Assume that $X_0\in \mH^p_1(TM)\cap L^\infty(TM)$ for some $p>1$ satisfies
$$
[\div X_0]^-\in L^\infty(M),
$$
and for each $k=1,\cdots, m$, $X_k\in C^2(TM)$.
Then there exists a unique $\nu$-almost everywhere stochastic
flows $\{x_t(x), x\in M\}_{t\geq 0}$ associated with SDE (\ref{SSDE})
in the sense of Definition \ref{Def1}.
\et
\begin{proof}
Let $X_{0,n}\in C^\infty(TM)$ be defined as in (\ref{Mol}).
Let $x_{n,t}(x)$ solve the following Stratonovich's SDE on $M$:
$$
\dif x_{n,t}(x)=X_{0,n}(x_{n,t}(x))\dif t+X_{k}(x_{n,t}(x))\circ \dif W^k_t,\ \ x_{n,0}=x.
$$
Then $x\mapsto x_{n,t}(x), t\geq 0$ defines a stochastic homeomorphism flow over $M$.
Moreover, by Proposition \ref{Pro1}
$$
(\nu\circ x_{n,t})(\dif x)=\beta_{n,t}(x)\nu(\dif x),
$$
where $\beta_{n,t}(x)$ satisfies by (\ref{Es22}) and (\ref{PP1}),  that for any $q\geq 1$,
\be
\sup_{n\in\mN}\sup_{(t,x)\in[0,T]\times M}\mE|\beta_{n,t}(x)|^q<+\infty.\label{Ep4}
\ee

Let us  set
$$
\Phi_{n,m}(x):=\sup_{t\in[0,T]}\dis^2(x_{n,t}(x),x_{m,t}(x))
$$
and
$$
A^\delta_{n,m}(x):=\log\left(\frac{\Phi_{n,m}(x)}{\delta}+1\right).
$$
If we choose
$$
\delta=\delta_{n,m}=\|X_{0,n}-X_{0,m}\|_1,
$$
then by Lemma \ref{Le2} and (\ref{Ep4}), we have
$$
\sup_{n,m}\mE \int_M A^{\delta_{n,m}}_{n,m}(x)\nu(\dif x)\leq C_0.
$$
Thus, by Chebyshev's inequality,  we have for any $R>0$,
\ce
\mE\int_M\Phi_{n,m}(x)\nu(\dif x)&=&\mE\int_M\Phi_{n,m}(x)\cdot
1_{\{A^{\delta_{n,m}}_{n,m}(x)>R\}}\nu(\dif x)
+\mE\int_M\Phi_{n,m}(x)\cdot 1_{\{A^{\delta_{n,m}}_{n,m}(x)\leq R\}}\nu(\dif x)\\
&\leq&\frac{(\mbox{diam}(M))^2\cdot C_0}{R}+\delta_{n,m}\cdot (e^{R}-1)\cdot\nu(M)\\
&=&\frac{(\mbox{diam}(M))^2\cdot C_0}{R}+\|X_{0,n}-X_{0,m}\|_1\cdot (e^{R}-1)\cdot\nu(M),
\de
where diam$(M):=\sup_{x,y\in M}\dis(x,y)$. From this, by Proposition \ref{Pro2},
we then obtain that
$$
\lim_{n,m\to\infty}\mE\int_M\sup_{t\in[0,T]}\dis^2(x_{n,t}(x),x_{m,t}(x))\nu(\dif x)=
\lim_{n,m\to\infty}\mE\int_M\Phi_{n,m}(x)\nu(\dif x)=0.
$$
Hence, for $\nu$-almost all $x\in M$, there exists a continuous
($\cF_t$)-adapted process $x_t(x)$ such that
\be
\lim_{n\to\infty}\mE\int_M\sup_{t\in[0,T]}\dis^2(x_{n,t}(x),x_{t}(x))\nu(\dif x)=0.\label{Es3}
\ee
By Lemma \ref{Le4} and (\ref{Ep4}), one finds that $x_t(x)$ satisfies {\bf (A)} and {\bf (B)}
 of Definition \ref{Def1}.
The uniqueness is a direct consequence of Lemma \ref{Le2}.
\end{proof}

We are now in a position to give

{\it Proof of Theorem \ref{Main}:}

Following the proof of Theorem \ref{Me}, we only need to
check {\bf (C)} of Definition \ref{Def1}.
Fix a $T>0$ and let
$$
\rho_{n}:=\exp\left\{\int^T_0\div X_{0,n}(x_{n,s}(x))\dif s
+\int^T_0\div X_k(x_{n,s}(x))\circ\dif W^k_s\right\}.
$$
By (\ref{Es2}), we have for any $q\geq 1$,
\be
\sup_{n\in\mN}\sup_{x\in\mR^d}\mE|\rho_n(x)|^q<+\infty.\label{PP3}
\ee
In view of (\ref{Ep4}) and (\ref{Es3}), by Lemma \ref{Le4}, we have
\ce
&&\lim_{n\to\infty}\mE\int^T_0\!\!\!\int_M|\div X_{0,n}(x_{n,s}(x))-\div X_0(x_s(x))|
\nu(\dif x)\dif s=0,\\
&&\lim_{n\to\infty}\mE\int_M
\left|\int^T_0(\div X_k(x_{n,s}(x))-\div X_k(x_s(x)))\circ\dif W^k_s\right|\nu(\dif x)=0.
\de
So, there is a subsequence still denoted by $n$ such that for almost all $(\omega,x)$,
\be
\lim_{n\to\infty}\rho_{n}(\omega,x)=\rho_T(\omega,x),\label{PP4}
\ee
where $\rho_T(x)$ is defined by (\ref{Rho}). By (\ref{PP3}) and (\ref{PP4}), we further have
for any $q\geq 1$,
\be
\lim_{n\to\infty}\mE\int_M|\rho_{n}(x)-\rho_T(x)|^q\nu(\dif x)=0.\label{U1}
\ee

Now, let $y_n(x)$ solve the following SDE
$$
\dif y_{n,t}(x)=-X_{0,n}(y_{n,t}(x))\dif t+X_k(y_{n,t}(x))\circ\dif W^{T,k}_t,\ \ y_n|_{t=0}=x,
$$
where  $W^T_t:=W_{T-t}-W_T$.
As in the proof of Theorem \ref{Me}, there exists a continuous ($\cF_t$)-adapted
process $y_{t}(x)$ such that
\be
\lim_{n\to\infty}\mE\int_M\sup_{t\in[0,T]}\dis(y_{n,t}(x),y_t(x))^2\nu(\dif x)=0.\label{U2}
\ee
It is well known that
$$
x^{-1}_{n,T}(x)=y_{n,T}(x).
$$
Thus, for any $f,g\in C(M)$, we have
\be
\int_M f(y_{n,T}(\omega,x))\cdot g(x)\nu(\dif x)
=\int_M f(x)\cdot g(x_{n,T}(\omega,x))\cdot\rho_{n}(\omega,x)\nu(\dif x),
\ \ P-a.s.\label{L01}
\ee
Let $\sC$ be a countable and dense subset of $C(M)$.
By (\ref{Es3}), (\ref{U1}) and (\ref{U2}),  if necessary, extracting a subsequence and
then taking limits $n\to\infty$ in $L^1(\Omega)$ for both sides of (\ref{L01}), we get
that for all $f,g\in\sC\subset C(M)$ and $P$-almost all $\omega\in\Omega$,
\be
\int_M f(y_{T}(\omega,x))\cdot g(x)\nu(\dif x)
=\int_M f(x)\cdot g(x_{T}(\omega,x))\cdot\rho_{T}(\omega,x)\nu(\dif x).\label{L1}
\ee
Since $\sC$ is countable, one may find a common null set $\Omega'\subset\Omega$
such that (\ref{L1}) holds for all $\omega\notin\Omega'$ and $f,g\in\sC$.
Thus, by Lemma \ref{Le33}, one sees that {\bf (C)} of Definition \ref{Def1} holds.


\begin{thebibliography}{999}

\bibitem{Am}Ambrosio, L.: Transport equation and
Cauchy problem for $BV$ vector fields.  Invent. Math.,  158  (2004),  no. 2, 227--260.

\bibitem{Ar}Arnold, L.: Random dynamical systems. Springer Monographs
in Mathematics. Springer-Verlag, Berlin, 1998.

\bibitem{B-C} Bishop, R.L.,  Crittenden, R.J.:   Geometry of manifolds. Academic Press Inc.,
New York,  1964.

\bibitem{Cr-De-Le}Crippa G. and De Lellis C.: Estimates and regularity results for
the DiPerna-Lions flow. J. reine angew. Math. 616 (2008), 15-46.

\bibitem{Di-Li}DiPerna R.J. and Lions P.L.: Ordinary differential equations, transport theory
and Sobolev spaces.  Invent. Math., 98,511-547(1989).

\bibitem{Du-Go-Lo}Dumas, H.S., Golse, F. and Lochak, P.: Multiphase averaging for generalized
flows on manifolds. Ergodic Theory and Dynamical System,  Vol. 53-63, 1994.

\bibitem{El}Elworthy, K.D.:  Stochastic differential equations on manifold.
London Math. Soc. Lecture Notes in Math., 70, Cambridge, University Press.

\bibitem{iw}Ikeda, N., Watanabe, S.:  Stochastic differential
equations and diffusion processes. North-Holland/Kodanska 1989.

\bibitem{K-N}Kobayashi, S., Nomizu K.: Differential geometry.
Interscience Publishers, Wiley, New York, 1963.

\bibitem{Ku}Kunita, H.: Stochastic flows and stochastic differential equations.
Cambridge, Cambridge University Press, 1990.

\bibitem{Le-Li}Le Bris, C.  and Lions, P.L. : Renormalized solutions of some transport equations
with partially $W^{1,1}$ velocities and applications.
Annali di Matematica, 183, 97-130(2004).

\bibitem{lee}Lee, H. M.:  Riemannian manifolds. GTM 176, Springer-Verlag
New York, 1997.

\bibitem{Re-Zh}Ren, J. and Zhang, X.: Limit theorems and large deviations
for stochastic flows on manifolds. Preprint.

\bibitem{St}Stein, E. M.: Singular integrals and differentiability properties of functions.
Princeton, N.J.,  Princeton University Press,  1970.

\bibitem{Str}  Stroock, D. W.:   An Introduction to the analysis of paths
on a Riemannian manifold. Math.\ Surveys and Monographs, 74 (AMS, 2000).

\bibitem{Zh2}Zhang, X.: Stochastic flows of SDEs
with irregular coefficients and stochastic transport equations. Bull. Sci. Math. France, Vol. 134, (2010)340-378.

\end{thebibliography}
\end{document}